\documentclass{article}
\usepackage[preprint]{neurips_2025} 

\usepackage[utf8]{inputenc} 
\usepackage[T1]{fontenc}    
\usepackage{hyperref}       
\usepackage{url}            
\usepackage{booktabs}       
\usepackage{amsfonts}       
\usepackage{nicefrac}       
\usepackage{microtype}      
\usepackage{xcolor}         
\usepackage[ruled,vlined,algosection,linesnumbered]{algorithm2e}

\usepackage{graphicx}
\usepackage{subfigure}

\usepackage{amsmath}
\usepackage{amssymb}
\usepackage{mathtools}
\usepackage{amsthm}
\theoremstyle{plain}

\theoremstyle{definition}

\theoremstyle{remark}

\usepackage[textsize=tiny]{todonotes}

\usepackage{multirow}
\input{helpers}
\title{High-dimensional Optimization with Low Rank Tensor Sampling and Local Search}
\author{%
    Konstantin Sozykin$^*$\\
    Skolkovo Institute\\of Science and Technology\\
    \texttt{konstantin.sozykin@skoltech.ru} \\
    \And
    Andrei Chertkov$^*$\\
    Artificial Intelligence Research
    Institute (AIRI),\\
    Skolkovo Institute of Science and Technology\\
    \texttt{a.chertkov@skoltech.ru}\\
    \And
    Anh-Huy Phan\\
    Skolkovo Institute\\of Science and Technology\\
    \texttt{a.phan@skoltech.ru} \\
    \And
    Ivan Oseledets\\
    Artificial Intelligence Research
    Institute (AIRI),\\
    Skolkovo Institute of Science and Technology\\
    \texttt{i.oseledets@skoltech.ru} \\
    \And
    Gleb Ryzhakov$^\dagger$\\
    Central University, Moscow, Russia\\
    \texttt{g.ryzhakov@centraluniversity.ru} \\
}
\begin{document}
\maketitle
\def\thefootnote{*}\footnotetext{Equal contribution}
\def\thefootnote{$\dagger$}\footnotetext{Corresponding author}
\def\thefootnote{\arabic{footnote}}
\begin{abstract}
We present a novel method called TESALOCS (TEnsor SAmpling and LOCal Search) for multidimensional optimization, combining the strengths of gradient-free discrete methods and gradient-based approaches. The discrete optimization in our method is based on low-rank tensor techniques, which, thanks to their low-parameter representation, enable efficient optimization of high-dimensional problems. For the second part, i.e., local search, any effective gradient-based method can be used, whether existing (such as quasi-Newton methods) or any other developed in the future. Our approach addresses the limitations of gradient-based methods, such as getting stuck in local optima; the limitations of discrete methods, which cannot be directly applied to continuous functions; and limitations of gradient-free methods that require large computational budgets. Note that we are not limited to a single type of low-rank tensor decomposition for discrete optimization, but for illustrative purposes, we consider a specific efficient low-rank tensor train decomposition. For 20 challenging 100-dimensional functions, we demonstrate that our method can significantly outperform results obtained with gradient-based methods like Conjugate Gradient, BFGS, SLSQP, and other methods, improving them by orders of magnitude with the same computing budget.
\end{abstract}
\section{Introduction}
    \label{s:intro}
    A persistent challenge in numerical optimization arises from the susceptibility of gradient-based methods to converge to local optima, particularly in high-dimensional search spaces where solution quality exhibits strong dependence on initialization.
While these methods demonstrate remarkable efficiency in local refinement, their global performance remains fundamentally constrained by the quality of starting points.
To overcome this limitation, we propose a hybrid paradigm that employs gradient-free strategies for identifying promising initial regions in the parameter space followed by gradient-based local optimization.

Contemporary approaches in this vein often utilize sampling mechanisms based on multivariate normal distributions~\cite{cmaes2011, NIPS2016_cmaes} or their mixture variants~\cite{9244595}, which theoretically enable systematic exploration of complex landscapes.
However, despite their mathematical elegance, such methods frequently underperform in practical high-dimensional scenarios due to the insufficient expressiveness of conventional sampling distributions to capture complex parameter correlations.

In the case of discrete parameter spaces, low rank tensor methods provide a powerful tool for building expressive yet low-parametric representations of high-dimensional models.
State-of-the-art algorithms like TTopt~\cite{sozykin2022ttopt} and PROTES~\cite{batsheva2023protes} exemplify this paradigm, leveraging the well-known Tensor Train (TT) format~\cite{oseledets2011tensor, cichocki2016tensor, cichocki2017tensor} to solve complex discrete black-box optimization problems ranging from quantum computing~\cite{paradezhenko2024probabilistic} to reinforcement learning~\cite{shetty2022tensor}.
The efficacy of these methods stems from the intrinsic properties of tensor networks like TT, Canonical Polyadic (CP)~\cite{harshman1970foundations}, and Hierarchical Tucker (HT)~\cite{hackbusch2009new}, which enable compact representations of high-dimensional data structures while maintaining computational tractability.
This combination of expressiveness and efficiency makes tensor networks particularly suitable for modern optimization problems where dimensionality scaling and resource limitations are fundamental constraints.

In this work, we establish a bridge between discrete low-rank tensor-based optimization and gradient-driven approach.
The motivation for our method is that conventional gradient methods are highly dependent on the starting point (approximation).
We solve this problem by iteratively sampling starting points from a discrete low-rank representation (surrogate model) of the optimized function that probabilistically encodes promising regions of the optimization landscape.
Moreover, the use of a low-rank approach allows us to compactly store the model, efficiently sample candidates for starting points, and dynamically update the model, since memory consumption and computational complexity grows linearly with dimension.
While we focus on TT-decomposition due to its numerical stability, it is important to note that the choice of a specific type of low-rank tensor decomposition is not fundamental, and other approaches, such as CP or HT, can also be successfully applied.

We summarize our main contributions as follows:
\begin{itemize}
    \item 
        We develop a new method TESALOCS for optimization (finding the minimum or maximum\footnote{
            Further, for concreteness, we will consider the minimization problem in this paper, while the proposed method can be applied to the discrete maximization problem without any modifications.
        } value) of multivariable functions based on a sampling approach from a dynamically updated probability distribution in the TT-format, combined with external optimization procedures such as the well-known second-order BFGS~\cite{Nocedal2006Numerical} or more heuristic zero-order methods like particle swarm optimization (PSA)~\cite{kennedy1995particle}.        
    \item
        We demonstrate superiority of TESALOCS over $6$ modern gradient-based methods and $3$ popular gradient-free approaches through extensive benchmarks on $20$ non-convex high-dimensional functions. 
        Using fixed hyperparameters across all experiments, our method achieves consistent performance improvements under equivalent computational budgets.
\end{itemize}
\section{Method}
    \label{sec:method}
    \begin{figure}[t!]
    \centering
    \includegraphics
        [width=0.8\linewidth]
        {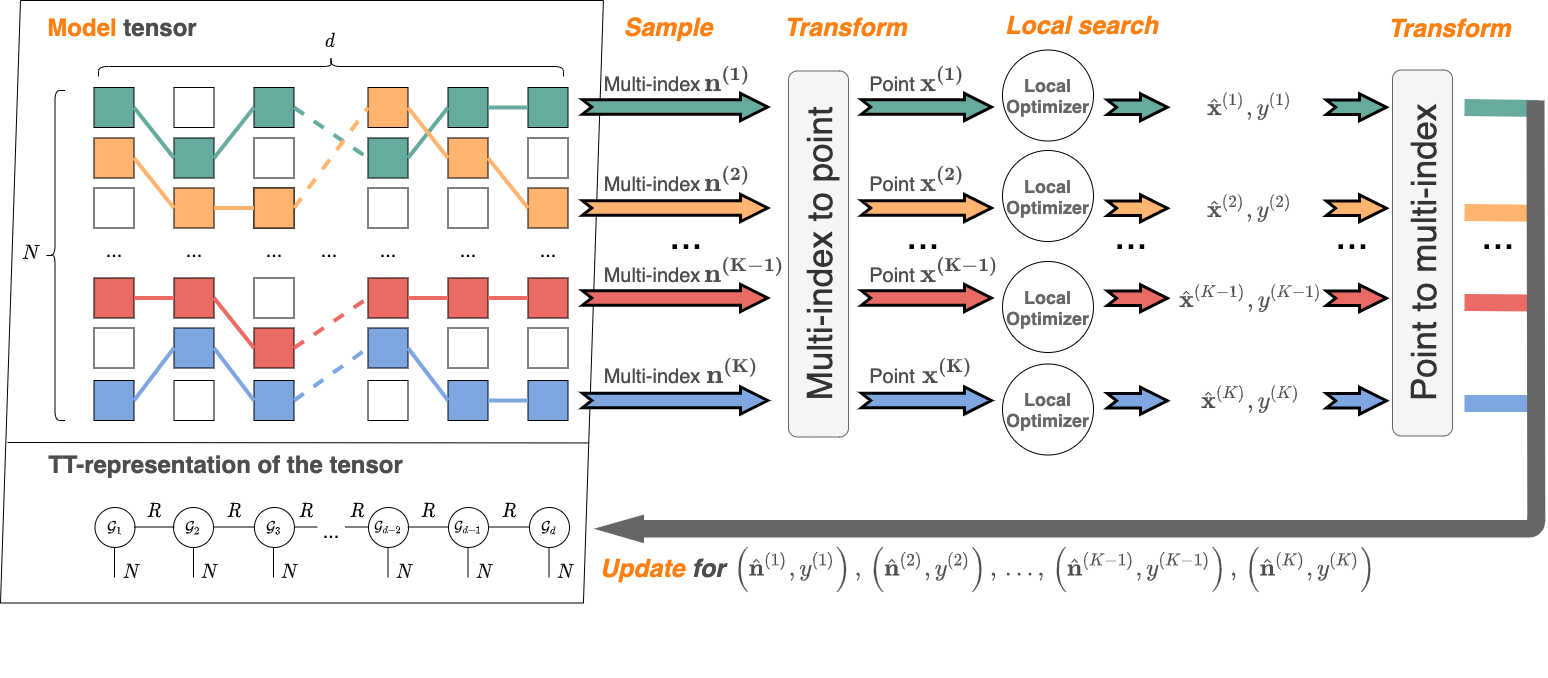}
    \vspace{-0.3cm}
    \caption{
        Schematic representation of the proposed optimization method TESALOCS.
    }
    \label{fig:tesalocs}
\vspace{-0.2cm}
\end{figure}

The TESALOCS algorithm is a novel method that decomposes high-dimensional optimization problem into two synergistic phases: global exploration through adaptive tensor sampling and local refinement via gradient-based or heuristic optimization.
The whole procedure is schematically illustrated in Figure~\ref{fig:tesalocs} and is formally described in Algorithm~\ref{alg:tesalocs}.
At its core, our method leverages low-rank tensor decompositions to construct a probabilistic discrete surrogate model that identifies promising regions in the parameter space while maintaining computational tractability. 
The algorithm operates by iteratively sampling candidate solutions from the TT-decomposition, refining them using a local search procedure, and updating the low-rank model in the TT-format based on the best-performing candidates.
This hybrid approach ensures both global exploration and local exploitation, making it highly effective for complex optimization problems.

Let us consider in more detail the mechanism of the TESALOCS method for optimization of $d$-dimensional function $f:\mathbb R^d \mapsto  \mathbb{R}$ in accordance with the main points indicated in Figure~\ref{fig:tesalocs}.

\paragraph{Model.}
We use a discrete surrogate model $\tens{T}$ in the TT-format to store information about the function $f(\mathbf{x})$ being optimized.
A $d$-dimensional tensor $\tens{T} \in \set{R}^{N_1 \times N_2 \times \ldots \times N_d}$ is expressed in the TT-format~\cite{oseledets2011tensor} if its elements can be represented as:
{
\begin{align} 
\label{eq:tt_repr}
\tens{T} [n_1, n_2, \ldots, n_d] = 
\sum_{r_1=1}^{R_1}
\sum_{r_2=1}^{R_2}
\cdots
\sum_{r_{d-1}=1}^{R_{d-1}}
    \tens{G}_1 [1, n_1, r_1]
    \tens{G}_2 [r_1, n_2, r_2]
    \ldots
    \tens{G}_d [r_{d-1}, n_d, 1],
\end{align}
}
where $(n_1, n_2, \ldots, n_d)$ is a multi-index (with $n_i = 1, 2, \ldots, N_i$ for $i = 1, 2, \ldots, d$), and the integers $R_{1}, R_{2}, \ldots, R_{d-1}$ (with the convention $R_{0} = R_{d} = 1$) controlling the expressiveness of the model and are called TT-ranks.
The three-dimensional tensors $\tens{G}_i \in \mathbb{R}^{R_{i-1} \times N_i \times R_i}$ (for $i = 1, 2, \ldots, d$) are known as TT-cores.
The TT-decomposition allows to represent a tensor in a compact and descriptive low-parameter form, which is linear in dimension~$d$, i.\,e., it has less than $d \cdot \max_{i=1,\ldots,d}(N_i R_i^2)$ parameters, hence the TT-format is recognized for its ability to reduce memory consumption and enable efficient linear algebra operations~\cite{cichocki2016tensor, cichocki2017tensor}.
Note that we start with a random non-negative tensor $\tens{T}$ with fixed rank $R_{1} = R_{2} = \ldots = R_{d-1} = R$, which is generated by assigning random non-negative elements to all its TT cores.
For simplicity, we use the equal number of discretization nodes for each dimension, i.e., $N_{1} = N_{2} = \ldots = N_{d} = N$.

\paragraph{Sample.}
Our algorithm proceeds iteratively within a given budget $M$, and in each iteration, a batch of $k$ candidate solutions is generated by sampling from the current tensor $\tens{T}$, which models a discrete probability distribution in the TT-format (that is, it corresponds to the probability for a given discrete point to be a global optimum).
The sampling is performed using an efficient procedure that has linear complexity in dimension $d$ (see~\cite{dolgov2020approximation} for details).

\paragraph{Transform.}
The generated discrete candidates $\mathbf{n}_{batch} = (\bf{n}^{(1)}, \bf{n}^{(2)}, \ldots, \bf{n}^{(k)})$ are projected onto a continuous grid using a projection operator $\textit{Pr}$ (``multi-index to point'' in Fig. \ref{fig:tesalocs}) and search space $\mathcal{X}$ defined as a tuple of $\mathbf{a}, \mathbf{b}, \mathbf{N}$ that represents lower and upper grid bound and a number of discretization nodes (``sample rate'') for search space.
In the case of uniform grid the projection operations $\textit{Pr}$ and $\textit{Pr}^{-1}$ (``point to multi-index'' in Fig. \ref{fig:tesalocs}) can be represented as follows:
\begin{align*}
    \textit{Pr}(\mathbf{\bf{n}^{(l)}, a, b, N}) =
        \frac
            {\bf{n}^{(l)}}
            {\mathbf{N} - 1}
            \cdot
            (\mathbf{b} - \mathbf{a}) + \mathbf{a},
    \quad
    \textit{Pr}^{-1}({\bf{x}^{(l)}}, \mathbf{a, b, N}) =    
        \frac
            {\bf{x}^{(l)} - \mathbf{a}}
            {\mathbf{b} - \mathbf{a}}
            \cdot
            (\mathbf{N} - 1).
\end{align*}
However, for different purposes, we may use a non-uniform grid or even a mix of different grids on different components of the input vector.
Note that we always have one-to-one mapping between multi-indices and corresponding float points.

\begin{algorithm}[t!]
\caption{TEnsor SAmpling and LOCal Search}
\label{alg:tesalocs}
\SetAlgoLined
\KwIn{
    Function $f(\mathbf{x})$, search space $\mathcal{X}$, budget $M$, discrete grid size $N$, TT-rank $r$, batch size $k$, elite size $k_{top}$, projectors $\textit{Pr}$ and $\textit{Pr}^{-1}$, local search algorithm \texttt{Loc}, SGD algorithm.
}
\KwOut{Optimal solution $\mathbf{x}^*$.}

Initialize cores $\{\mathcal G_i\}_{i=1}^d$ of the TT-decomposition~$\mathcal T$ with all $N_{1}, N_{2}, \ldots, N_{d}$ equal to $N$ and all $R_{1}, R_{2}, \ldots, R_{d-1}$ equal to $r$ \;

Set $m = 0$ and $\mathcal P =\emptyset$\;

\While{$m < M$}{
    Generate discrete candidates $\mathbf{n}_{batch} = \{\mathbf{n}^{(l)}\}_{l=1}^k$ via sampling from $\mathcal T$ in the TT-format\;
    
    Project discrete candidates on the continuous grid: $\mathbf{x}_{batch} = \textit{Pr} (\mathbf{n}_{batch},\mathcal{X}) $\;
    
    Make local search with the given algorithm, using  $\mathbf{x}_{batch}$ as initials:
    $\{\widehat{\mathbf{x}}_{batch}, y_{batch}, m_{\text{loc}}\}=\texttt{Loc}(f, \mathbf{x}_{batch})$\;

    Select the best $k_{top}$ candidates: $\mathbf{x}_{\text{top}}=\text{Top}(\widehat{\mathbf{x}}_{batch}, y_{batch}, k_{top})$\;
    
    Merge populations: $\mathcal{P} \leftarrow \mathcal{P} \cup \mathbf{x}_{\text{top}}$\;
    
    Unproject back: $\mathbf{n}_{\text{top}}=\textit{Pr}^{-1}(\mathbf{x}_{\text{top}}, \mathcal{X})$\;

    Make SGD step to update the cores using Log Likelihood loss $\mathcal{L}$ on the obtained indices: 
    $\{\mathcal G_i\}_{i=1}^d\gets SGD(\{\mathcal G_i\}_{i=1}^d, \mathbf{n}_{\text{top}}) $\;
    
    Update current budget: $m \leftarrow m + m_{\text{loc}}$\;
}

$\mathbf{x}^* = \arg\max_{\mathbf{x} \in \mathcal{P}} f(\mathbf{x})$\;


\end{algorithm}

\paragraph{Local search.}
A local search algorithm is applied to refine the projected candidates $\mathbf{x}_{batch} = (\bf{x}^{(1)}, \bf{x}^{(2)}, \ldots, \bf{x}^{(k)})$ to exploit high-potential regions identified by the global tensor sampler.
We can employ any gradient-based method (e.g., BFGS), or even derivative-free heuristics (e.g., PSO) in the case of non-differentiable or noisy objectives.
Each candidate $\bf{x}^{(l)}$ ($l = 1, 2, \ldots, k$) serves as an initial point for an independent local optimization run.
Based on the results of the local optimizer's work, we obtain a batch of refined continuous points $\widehat{\mathbf{x}}_{batch}$, the corresponding values of the objective function at these points $y_{batch}$, and the actual budget spent on this launch $m_{loc}$.
Note that these calculations can be performed in parallel if necessary to maximize computational efficiency.

\paragraph{Update.}
Following local refinement, the top-\( k \) solutions are selected from the batch $\mathbf{x}_{batch}$ based on their objective values, forming the elite subset $\mathbf{x}_{top}$.
Then we project them back into discrete space and use selected discrete $k_{top}$ candidates $\mathbf{n}_{top}$ in SGD with log-likelihood objective, which can be easily and efficiency computed using equation~\eqref{eq:tt_repr}:
\begin{align*}
    \mathcal{L} = -\sum_{\mathbf{n}\in\mathbf{n}_{top}} \log p_{\mathcal{T}}(\mathbf{n}),
\end{align*}
where $p_{\mathcal{T}}(\mathbf{n})$ is the probability density encoded in the TT-cores (i.e. normalized value of the TT-tensor). By doing gradient ascend with such loss we will increase likelihood of points $\mathbf{n}_{top}$ refined with local optimization.
Note that such update mechanism focuses the tensor sampler on promising subspaces while preserving rank constraints to maintain linear memory scaling.
\section{Numerical experiments}
    \label{sec:experiments}
    \begin{figure}[t!]
    \centering
    \subfigure{
        \centering
        \includegraphics[width=0.31\linewidth]{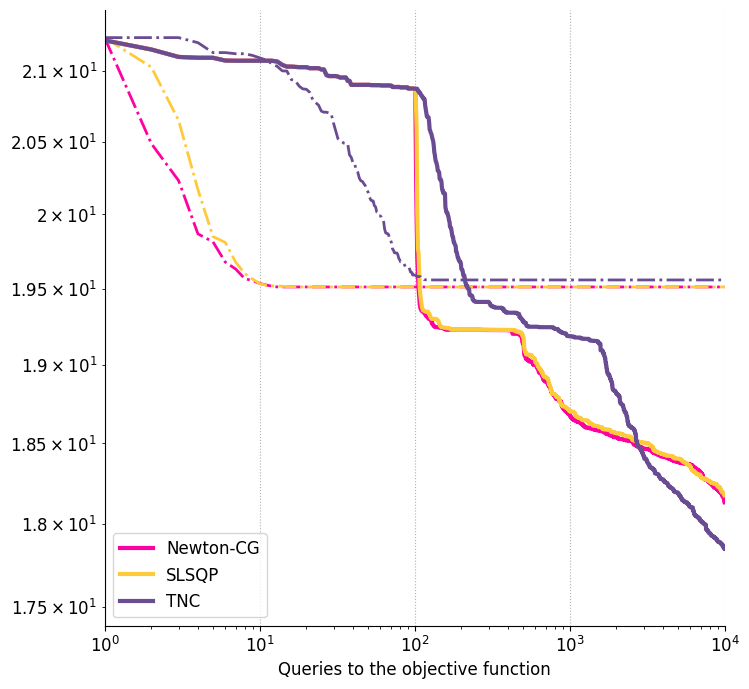}
    }
    \subfigure{
        \centering
        \includegraphics[width=0.31\linewidth]{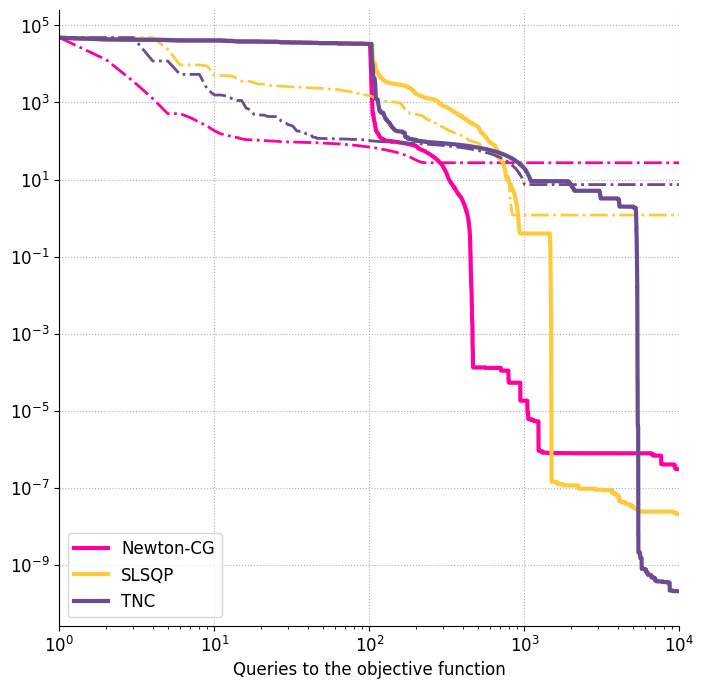}
    }
    \subfigure{
        \centering
        \includegraphics[width=0.31\linewidth]{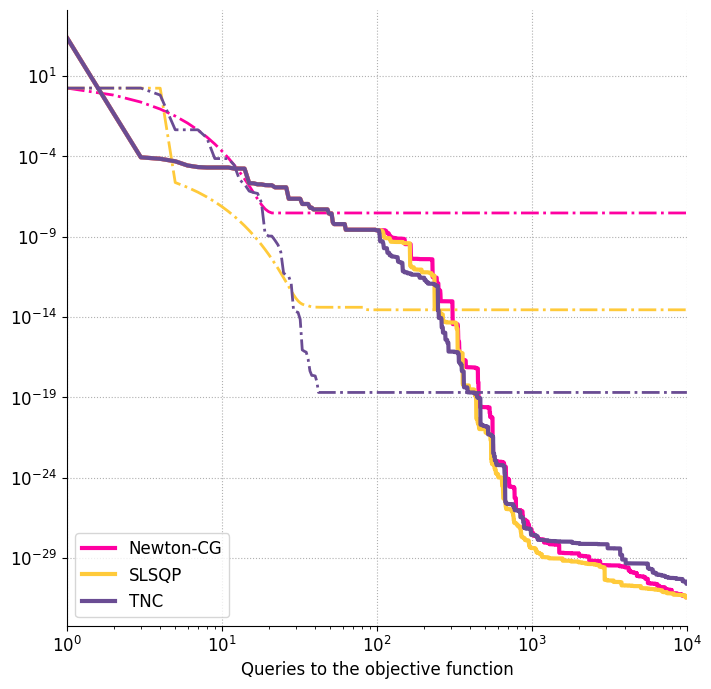}
    }
    \renewcommand{\thesubfigure}{(a)}
    \subfigure[Ackley function]{
        \centering
        \includegraphics[width=0.31\linewidth]{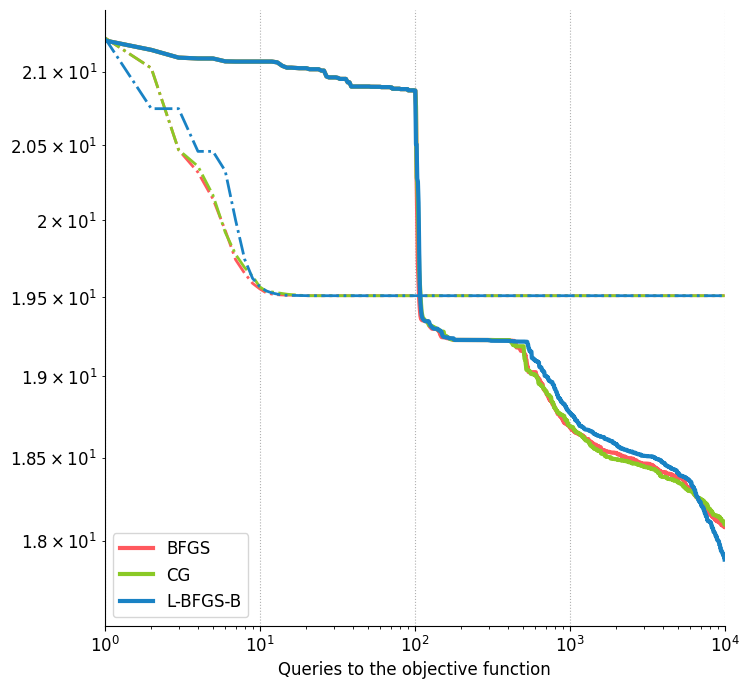}
    }
    \renewcommand{\thesubfigure}{(b)}
    \subfigure[Rosenbrock function]{
        \centering
        \includegraphics[width=0.31\linewidth]{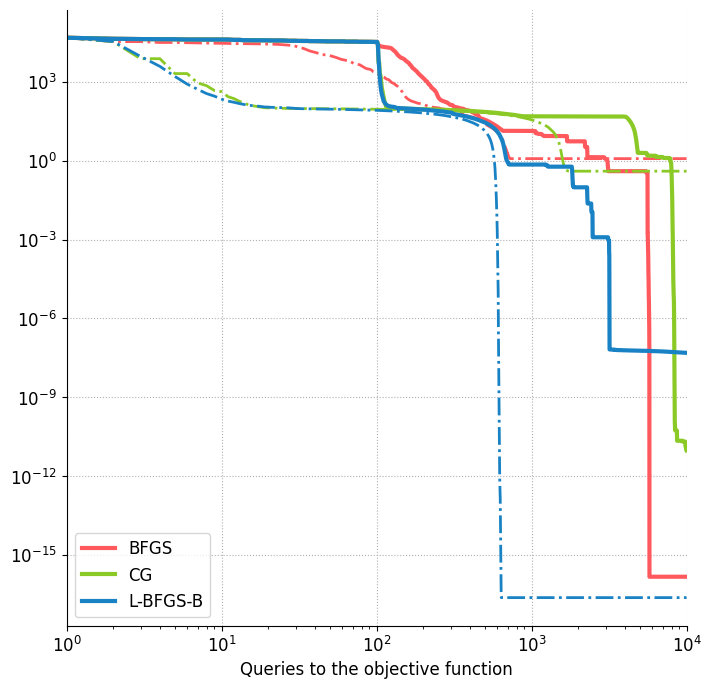}
    }
    \renewcommand{\thesubfigure}{(c)}
    \subfigure[Yang function]{
        \centering
        \includegraphics[width=0.31\linewidth]{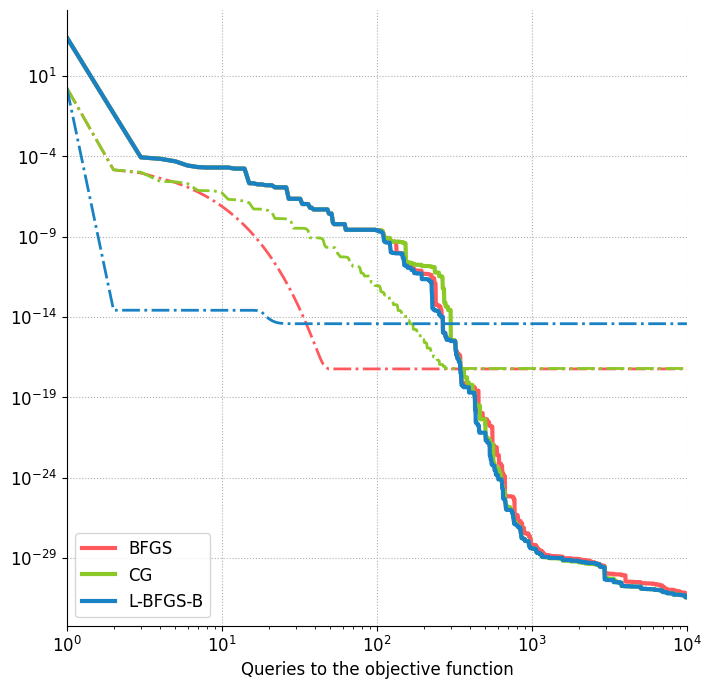}
    }
    \caption{
        Comparative minimization results for six optimization algorithms (Newton-CG, SLSQP, TNC, BFGS, CG, and L-BFGS-B) applied to the $100$-dimensional Ackley (left), Rosenbrock (middle), and Yang (right) functions.
        For each method, the graphs display the average minimum value over $10$ runs, depending on the number of queries made to the target function.        
        The conventional baseline approach (random initialization without enhancements) is indicated by dashed lines, while results from our proposed TESALOCS method are shown with solid lines. 
    }
    \label{fig:results_conv}
\vspace{-0.5cm}
\end{figure}

\begin{table}
\caption{Optimization results for gradient-based methods BFGS, CG, and L-BFGS-B.}
\label{tab:experiments_func1}
\centering
\begin{tabular}{llllllll}
\toprule
\multirow{2}{*}{Function}
&
& \multicolumn{2}{c}{BFGS}
& \multicolumn{2}{c}{CG}
& \multicolumn{2}{c}{L-BFGS-B}
\\
\cmidrule(r){3-4}
\cmidrule(r){5-6}
\cmidrule(r){7-8}
&
& Random & TESALOCS
& Random & TESALOCS
& Random & TESALOCS
\\
\midrule
{\bf\# of best results}
&
&  9 & {\bf19}
&  7 & {\bf19}
& 10 & {\bf18}
\\
\midrule
\multirow{2}{*}{Ackley}
& $E$
& 2.0e+01 & {\bf1.8e+01}
& 2.0e+01 & {\bf1.8e+01}
& 2.0e+01 & {\bf1.8e+01}
\\
& $\sigma$
& 8.2e-02 & 5.1e-01
& 8.2e-02 & 4.2e-01
& 8.2e-02 & 5.3e-01
\\
\midrule
\multirow{2}{*}{Alpine}
& $E$
& 3.8e-04 & {\bf1.8e-09}
& 2.2e+01 & {\bf7.5e+00}
& 6.2e+00 & {\bf8.0e-01}
\\
& $\sigma$
& 7.8e-04 & 1.1e-09
& 9.2e+00 & 3.2e+00
& 2.9e+00 & 7.3e-01
\\
\midrule
\multirow{2}{*}{Chung}
& $E$
& {\bf1.8e-22} & {\bf1.2e-09}
& 6.9e+03 & {\bf1.7e-11}
& {\bf2.3e-17} & {\bf1.6e-10}
\\
& $\sigma$
& 2.9e+00 & 7.3e-01
& 8.2e+03 & 3.4e-11
& 8.2e+03 & 3.4e-11
\\
\midrule
\multirow{2}{*}{Dixon}
& $E$
& {\bf6.7e-01} & {\bf6.7e-01}
& {\bf6.7e-01} & {\bf6.7e-01}
& {\bf6.7e-01} & {\bf6.7e-01}
\\
& $\sigma$
& 3.6e-15 & 0.0e+00
& 1.1e-15 & 0.0e+00
& 2.7e-16 & 0.0e+00
\\
\midrule
\multirow{2}{*}{Exp}
& $E$
& {\bf0.0e+00} & {\bf0.0e+00}
& {\bf0.0e+00} & {\bf0.0e+00}
& {\bf0.0e+00} & {\bf0.0e+00}
\\
& $\sigma$
& 2.7e-16 & 0.0e+00
& 2.7e-16 & 0.0e+00
& 2.7e-16 & 0.0e+00
\\
\midrule
\multirow{2}{*}{Griewank}
& $E$
& {\bf0.0e+00} & {\bf0.0e+00}
& 7.4e-04 & {\bf2.7e-10}
& {\bf0.0e+00} & {\bf0.0e+00}
\\
& $\sigma$
& 2.7e-16 & 0.0e+00
& 2.2e-03 & 2.4e-10
& 2.2e-03 & 2.4e-10
\\
\midrule
\multirow{2}{*}{Pathological}
& $E$
& 3.8e+01 & {\bf3.0e+01}
& 4.0e+01 & {\bf3.2e+01}
& 3.9e+01 & {\bf3.2e+01}
\\
& $\sigma$
& 1.4e+00 & 1.5e+00
& 1.3e+00 & 2.0e+00
& 1.6e+00 & 1.5e+00
\\
\midrule
\multirow{2}{*}{Pinter}
& $E$
& 1.5e+05 & {\bf6.0e+04}
& 1.7e+05 & {\bf6.6e+04}
& 1.7e+05 & {\bf6.4e+04}
\\
& $\sigma$
& 2.1e+04 & 3.7e+03
& 2.0e+04 & 9.2e+03
& 2.0e+04 & 6.8e+03
\\
\midrule
\multirow{2}{*}{Powell}
& $E$
& {\bf1.3e-18} & 6.2e-08
& {\bf1.2e-17} & 2.3e-08
& {\bf2.6e-11} & 2.8e-08
\\
& $\sigma$
& 1.2e-18 & 3.8e-08
& 1.0e-17 & 2.0e-08
& 2.0e-11 & 1.5e-08
\\
\midrule
\multirow{2}{*}{Qing}
& $E$
& {\bf8.9e-23} & {\bf3.8e-14}
& {\bf3.8e-26} & {\bf2.7e-13}
& {\bf1.1e-16} & {\bf1.6e-09}
\\
& $\sigma$
& 2.0e-11 & 1.5e-08
& 2.0e-11 & 1.5e-08
& 2.0e-11 & 1.5e-08
\\
\midrule
\multirow{2}{*}{Rastrigin}
& $E$
& 8.4e+02 & {\bf3.3e+02}
& 8.5e+02 & {\bf3.2e+02}
& 8.5e+02 & {\bf2.9e+02}
\\
& $\sigma$
& 7.9e+01 & 4.5e+01
& 7.4e+01 & 5.0e+01
& 7.4e+01 & 4.6e+01
\\
\midrule
\multirow{2}{*}{Rosenbrock}
& $E$
& 1.2e+00 & {\bf1.5e-16}
& 4.0e-01 & {\bf1.0e-11}
& {\bf4.8e-17} & 2.3e-08
\\
& $\sigma$
& 1.8e+00 & 4.5e-16
& 1.2e+00 & 1.8e-11
& 3.9e-17 & 4.7e-08
\\
\midrule
\multirow{2}{*}{Salomon}
& $E$
& 5.7e+01 & {\bf1.5e+01}
& 5.7e+01 & {\bf1.0e+01}
& 4.6e+01 & {\bf9.0e-15}
\\
& $\sigma$
& 2.5e+00 & 1.0e+01
& 2.5e+00 & 9.2e+00
& 2.2e+01 & 4.6e-15
\\
\midrule
\multirow{2}{*}{Schaffer}
& $E$
& 4.5e+01 & {\bf4.0e+01}
& 4.5e+01 & {\bf4.0e+01}
& 4.5e+01 & {\bf4.0e+01}
\\
& $\sigma$
& 7.3e-01 & 1.4e+00
& 7.2e-01 & 1.4e+00
& 7.2e-01 & 1.4e+00
\\
\midrule
\multirow{2}{*}{Sphere}
& $E$
& {\bf7.6e-34} & {\bf3.1e-34}
& {\bf6.4e-39} & {\bf8.5e-19}
& {\bf8.9e-59} & {\bf5.8e-37}
\\
& $\sigma$
& 7.2e-01 & 1.4e+00
& 7.2e-01 & 1.4e+00
& 7.2e-01 & 1.4e+00
\\
\midrule
\multirow{2}{*}{Squares}
& $E$
& {\bf1.8e-35} & {\bf1.3e-20}
& {\bf1.4e-33} & {\bf5.3e-12}
& {\bf2.7e-16} & {\bf4.2e-10}
\\
& $\sigma$
& 7.2e-01 & 1.4e+00
& 7.2e-01 & 1.4e+00
& 7.2e-01 & 1.4e+00
\\
\midrule
\multirow{2}{*}{Trid}
& $E$
& 4.2e-08 & {\bf0.0e+00}
& 1.3e-05 & {\bf0.0e+00}
& 1.5e-06 & {\bf0.0e+00}
\\
& $\sigma$
& 3.8e-08 & 0.0e+00
& 6.0e-06 & 0.0e+00
& 1.0e-06 & 0.0e+00
\\
\midrule
\multirow{2}{*}{Trigonometric}
& $E$
& 3.9e+03 & {\bf3.0e+01}
& 1.3e+01 & {\bf7.4e-01}
& 2.2e+00 & {\bf2.4e-01}
\\
& $\sigma$
& 8.1e+02 & 9.1e+01
& 1.3e+01 & 9.8e-01
& 5.9e+00 & 7.1e-01
\\
\midrule
\multirow{2}{*}{Wavy}
& $E$
& 5.9e-01 & {\bf3.6e-01}
& 5.9e-01 & {\bf3.7e-01}
& 5.9e-01 & {\bf3.6e-01}
\\
& $\sigma$
& 3.4e-02 & 1.8e-02
& 3.4e-02 & 2.4e-02
& 3.4e-02 & 2.8e-02
\\
\midrule
\multirow{2}{*}{Yang}
& $E$
& {\bf5.9e-18} & {\bf4.9e-32}
& {\bf6.2e-18} & {\bf4.2e-32}
& {\bf3.8e-15} & {\bf4.2e-32}
\\
& $\sigma$
& 3.4e-02 & 2.8e-02
& 3.4e-02 & 2.8e-02
& 3.4e-02 & 2.8e-02
\\
\bottomrule
\end{tabular}
\end{table}

\begin{table}
\caption{Optimization results for gradient-based methods Newton-CG, SLSQP, and TNC.}
\label{tab:experiments_func2}
\centering
\begin{tabular}{llllllll}
\toprule
\multirow{2}{*}{Function}
&
& \multicolumn{2}{c}{Newton-CG}
& \multicolumn{2}{c}{SLSQP}
& \multicolumn{2}{c}{TNC}
\\
\cmidrule(r){3-4}
\cmidrule(r){5-6}
\cmidrule(r){7-8}
&
& Random & TESALOCS
& Random & TESALOCS
& Random & TESALOCS
\\
\midrule
{\bf\# of best results}
&
&  6 & {\bf20}
&  9 & {\bf15}
&  8 & {\bf20}
\\
\midrule
\multirow{2}{*}{Ackley}
& $E$
& 2.0e+01 & {\bf1.8e+01}
& 2.0e+01 & {\bf1.8e+01}
& 2.0e+01 & {\bf1.8e+01}
\\
& $\sigma$
& 8.3e-02 & 3.8e-01
& 8.2e-02 & 3.4e-01
& 9.2e-02 & 7.9e-01
\\
\midrule
\multirow{2}{*}{Alpine}
& $E$
& 9.1e+01 & {\bf2.5e+01}
& {\bf6.9e-13} & 3.1e-05
& 4.7e+01 & {\bf4.9e-04}
\\
& $\sigma$
& 2.1e+01 & 5.8e+00
& 1.7e-13 & 5.1e-06
& 1.8e+01 & 7.8e-04
\\
\midrule
\multirow{2}{*}{Chung}
& $E$
& {\bf6.6e-10} & {\bf4.6e-11}
& {\bf6.9e-13} & 3.7e-08
& {\bf4.1e-17} & {\bf2.1e-11}
\\
& $\sigma$
& 1.8e+01 & 7.8e-04
& 1.3e-12 & 5.3e-08
& 1.3e-12 & 5.3e-08
\\
\midrule
\multirow{2}{*}{Dixon}
& $E$
& 7.6e+00 & {\bf6.7e-01}
& {\bf6.7e-01} & {\bf6.7e-01}
& {\bf6.7e-01} & {\bf6.7e-01}
\\
& $\sigma$
& 1.9e+01 & 0.0e+00
& 5.8e-16 & 1.1e-07
& 7.1e-16 & 0.0e+00
\\
\midrule
\multirow{2}{*}{Exp}
& $E$
& {\bf0.0e+00} & {\bf0.0e+00}
& 4.0e-01 & {\bf0.0e+00}
& {\bf2.0e-16} & {\bf0.0e+00}
\\
& $\sigma$
& 7.1e-16 & 0.0e+00
& 4.9e-01 & 0.0e+00
& 4.9e-01 & 0.0e+00
\\
\midrule
\multirow{2}{*}{Griewank}
& $E$
& {\bf3.8e-14} & {\bf8.6e-15}
& {\bf0.0e+00} & {\bf0.0e+00}
& 7.9e-03 & {\bf1.5e-03}
\\
& $\sigma$
& 4.9e-01 & 0.0e+00
& 4.9e-01 & 0.0e+00
& 9.8e-03 & 3.0e-03
\\
\midrule
\multirow{2}{*}{Pathological}
& $E$
& 4.2e+01 & {\bf3.2e+01}
& 4.1e+01 & {\bf3.2e+01}
& 3.9e+01 & {\bf3.1e+01}
\\
& $\sigma$
& 2.3e+00 & 1.8e+00
& 5.0e+00 & 1.7e+00
& 1.4e+00 & 1.6e+00
\\
\midrule
\multirow{2}{*}{Pinter}
& $E$
& 1.7e+05 & {\bf7.0e+04}
& 1.3e+05 & {\bf4.8e+04}
& 1.5e+05 & {\bf5.1e+04}
\\
& $\sigma$
& 2.0e+04 & 3.2e+03
& 1.4e+04 & 8.4e+03
& 2.2e+04 & 5.7e+03
\\
\midrule
\multirow{2}{*}{Powell}
& $E$
& 7.2e-08 & {\bf8.7e-09}
& {\bf6.2e-11} & 3.2e-06
& {\bf1.6e-14} & {\bf6.9e-13}
\\
& $\sigma$
& 5.3e-08 & 3.1e-09
& 4.7e-11 & 1.6e-06
& 4.7e-11 & 1.6e-06
\\
\midrule
\multirow{2}{*}{Qing}
& $E$
& {\bf1.4e-15} & {\bf5.1e-09}
& {\bf1.2e-14} & 3.0e-08
& {\bf1.1e-14} & {\bf5.5e-10}
\\
& $\sigma$
& 4.7e-11 & 1.6e-06
& 9.7e-15 & 1.9e-08
& 9.7e-15 & 1.9e-08
\\
\midrule
\multirow{2}{*}{Rastrigin}
& $E$
& 8.6e+02 & {\bf3.4e+02}
& 8.5e+02 & {\bf3.2e+02}
& 8.8e+02 & {\bf3.1e+02}
\\
& $\sigma$
& 8.1e+01 & 5.7e+01
& 7.3e+01 & 5.3e+01
& 8.8e+01 & 4.4e+01
\\
\midrule
\multirow{2}{*}{Rosenbrock}
& $E$
& 1.9e+01 & {\bf1.2e-07}
& 1.2e+00 & {\bf1.2e-08}
& 7.8e+00 & {\bf1.1e-10}
\\
& $\sigma$
& 3.6e+01 & 2.0e-07
& 1.8e+00 & 5.1e-09
& 4.9e+00 & 7.7e-11
\\
\midrule
\multirow{2}{*}{Salomon}
& $E$
& 5.7e+01 & {\bf2.1e+01}
& 5.7e+01 & {\bf2.5e+01}
& 5.6e+01 & {\bf1.4e+01}
\\
& $\sigma$
& 2.5e+00 & 8.6e+00
& 2.4e+00 & 2.1e+00
& 2.3e+00 & 1.1e+01
\\
\midrule
\multirow{2}{*}{Schaffer}
& $E$
& 4.6e+01 & {\bf4.0e+01}
& 4.5e+01 & {\bf4.0e+01}
& 4.5e+01 & {\bf3.9e+01}
\\
& $\sigma$
& 7.1e-01 & 9.0e-01
& 7.3e-01 & 1.5e+00
& 1.0e+00 & 1.5e+00
\\
\midrule
\multirow{2}{*}{Sphere}
& $E$
& {\bf8.6e-18} & {\bf2.9e-31}
& {\bf0.0e+00} & {\bf0.0e+00}
& {\bf1.7e-16} & {\bf1.0e-32}
\\
& $\sigma$
& 1.0e+00 & 1.5e+00
& 1.0e+00 & 1.5e+00
& 1.0e+00 & 1.5e+00
\\
\midrule
\multirow{2}{*}{Squares}
& $E$
& {\bf1.3e-19} & {\bf3.8e-11}
& {\bf1.0e-17} & 6.0e-08
& {\bf1.2e-15} & {\bf2.2e-10}
\\
& $\sigma$
& 1.0e+00 & 1.5e+00
& 5.0e-18 & 1.5e-08
& 5.0e-18 & 1.5e-08
\\
\midrule
\multirow{2}{*}{Trid}
& $E$
& 9.2e-05 & {\bf0.0e+00}
& 8.3e-08 & {\bf0.0e+00}
& 1.4e-07 & {\bf0.0e+00}
\\
& $\sigma$
& 1.1e-04 & 0.0e+00
& 6.1e-08 & 0.0e+00
& 8.9e-08 & 0.0e+00
\\
\midrule
\multirow{2}{*}{Trigonometric}
& $E$
& 1.1e+02 & {\bf1.5e+01}
& 1.1e+03 & {\bf1.3e+02}
& 1.0e+02 & {\bf3.3e-01}
\\
& $\sigma$
& 1.5e+02 & 1.8e+01
& 3.5e+02 & 1.1e+02
& 1.6e+02 & 6.7e-01
\\
\midrule
\multirow{2}{*}{Wavy}
& $E$
& 6.3e-01 & {\bf3.8e-01}
& 5.9e-01 & {\bf3.5e-01}
& 5.9e-01 & {\bf3.5e-01}
\\
& $\sigma$
& 6.6e-02 & 3.2e-02
& 3.4e-02 & 2.2e-02
& 3.5e-02 & 2.9e-02
\\
\midrule
\multirow{2}{*}{Yang}
& $E$
& 2.9e-08 & {\bf1.8e-32}
& {\bf1.5e-13} & {\bf4.2e-32}
& {\bf2.0e-19} & {\bf1.4e-31}
\\
& $\sigma$
& 1.1e-08 & 2.8e-32
& 1.1e-08 & 2.8e-32
& 1.1e-08 & 2.8e-32
\\
\bottomrule
\end{tabular}
\end{table}

\begin{table}
\caption{Optimization results for gradient-free methods NoisyBandit, PSO, and SPSA.}
\label{tab:experiments_func3}
\centering
\begin{tabular}{llllllll}
\toprule
\multirow{2}{*}{Function}
&
& \multicolumn{2}{c}{NoisyBandit}
& \multicolumn{2}{c}{PSO}
& \multicolumn{2}{c}{SPSA}
\\
\cmidrule(r){3-4}
\cmidrule(r){5-6}
\cmidrule(r){7-8}
&
& Random & TESALOCS
& Random & TESALOCS
& Random & TESALOCS
\\
\midrule
{\bf\# of best results}
&
&  1 & {\bf20}
&  3 & {\bf18}
&  0 & {\bf20}
\\
\midrule
\multirow{2}{*}{Ackley}
& $E$
& 1.8e+01 & {\bf1.4e+01}
& 1.9e+01 & {\bf1.4e+01}
& 2.1e+01 & {\bf1.8e+01}
\\
& $\sigma$
& 2.8e-01 & 1.1e+00
& 2.6e-01 & 1.1e+00
& 8.3e-02 & 1.4e+00
\\
\midrule
\multirow{2}{*}{Alpine}
& $E$
& 1.4e+02 & {\bf5.7e+01}
& 1.3e+02 & {\bf5.6e+01}
& 3.0e+02 & {\bf9.1e+01}
\\
& $\sigma$
& 4.4e+00 & 9.0e+00
& 8.6e+00 & 6.6e+00
& 1.3e+01 & 3.1e+01
\\
\midrule
\multirow{2}{*}{Chung}
& $E$
& 6.2e+05 & {\bf7.1e+04}
& 5.6e+05 & {\bf7.1e+04}
& 1.1e+07 & {\bf8.2e+05}
\\
& $\sigma$
& 1.2e+05 & 5.4e+04
& 1.2e+05 & 6.0e+04
& 1.8e+06 & 7.2e+05
\\
\midrule
\multirow{2}{*}{Dixon}
& $E$
& 3.1e+06 & {\bf6.4e+05}
& 3.2e+06 & {\bf7.4e+05}
& 3.8e+07 & {\bf5.6e+06}
\\
& $\sigma$
& 4.7e+05 & 1.5e+05
& 6.4e+05 & 3.0e+05
& 6.7e+06 & 4.1e+06
\\
\midrule
\multirow{2}{*}{Exp}
& $E$
& 9.8e-01 & {\bf7.1e-01}
& 9.8e-01 & {\bf7.1e-01}
& 1.0e+00 & {\bf9.6e-01}
\\
& $\sigma$
& 5.4e-03 & 1.0e-01
& 7.0e-03 & 6.2e-02
& 4.9e-07 & 2.7e-02
\\
\midrule
\multirow{2}{*}{Griewank}
& $E$
& 2.3e+01 & {\bf8.7e+00}
& 2.0e+01 & {\bf7.4e+00}
& 8.2e+01 & {\bf2.2e+01}
\\
& $\sigma$
& 1.8e+00 & 1.8e+00
& 1.5e+00 & 2.1e+00
& 7.0e+00 & 8.9e+00
\\
\midrule
\multirow{2}{*}{Pathological}
& $E$
& 4.8e+01 & {\bf4.2e+01}
& 4.3e+01 & {\bf3.8e+01}
& 5.0e+01 & {\bf4.2e+01}
\\
& $\sigma$
& 3.6e-01 & 1.3e+00
& 7.3e-01 & 1.4e+00
& 6.0e-01 & 1.8e+00
\\
\midrule
\multirow{2}{*}{Pinter}
& $E$
& 1.0e+05 & {\bf6.1e+04}
& 1.0e+05 & {\bf5.7e+04}
& 2.4e+05 & {\bf9.8e+04}
\\
& $\sigma$
& 3.3e+03 & 9.9e+03
& 3.6e+03 & 7.3e+03
& 2.1e+04 & 2.4e+04
\\
\midrule
\multirow{2}{*}{Powell}
& $E$
& 2.1e-02 & {\bf1.4e-03}
& 2.0e-02 & {\bf3.0e-04}
& 3.3e+00 & {\bf2.6e-02}
\\
& $\sigma$
& 1.4e-02 & 1.1e-03
& 5.9e-03 & 1.7e-04
& 1.0e+00 & 2.1e-02
\\
\midrule
\multirow{2}{*}{Qing}
& $E$
& 8.1e+11 & {\bf2.7e+10}
& 9.8e+11 & {\bf3.7e+10}
& 5.0e+20 & {\bf7.2e+10}
\\
& $\sigma$
& 3.3e+11 & 5.3e+10
& 7.1e+11 & 5.8e+10
& 0.0e+00 & 8.9e+10
\\
\midrule
\multirow{2}{*}{Rastrigin}
& $E$
& 1.1e+03 & {\bf5.4e+02}
& 1.1e+03 & {\bf4.6e+02}
& 1.8e+03 & {\bf7.6e+02}
\\
& $\sigma$
& 5.3e+01 & 1.2e+02
& 4.9e+01 & 7.3e+01
& 1.0e+02 & 2.1e+02
\\
\midrule
\multirow{2}{*}{Rosenbrock}
& $E$
& 6.2e+03 & {\bf2.1e+03}
& 7.0e+03 & {\bf1.6e+03}
& 4.8e+04 & {\bf7.8e+03}
\\
& $\sigma$
& 9.1e+02 & 7.4e+02
& 8.3e+02 & 5.2e+02
& 4.3e+03 & 4.3e+03
\\
\midrule
\multirow{2}{*}{Salomon}
& $E$
& 3.0e+01 & {\bf1.9e+01}
& 2.9e+01 & {\bf1.8e+01}
& 5.8e+01 & {\bf3.3e+01}
\\
& $\sigma$
& 1.7e+00 & 2.1e+00
& 1.2e+00 & 2.3e+00
& 2.6e+00 & 6.7e+00
\\
\midrule
\multirow{2}{*}{Schaffer}
& $E$
& 4.8e+01 & {\bf3.8e+01}
& {\bf3.2e+01} & 3.7e+01
& 5.0e+01 & {\bf4.2e+01}
\\
& $\sigma$
& 4.9e-01 & 1.5e+00
& 6.1e+00 & 1.3e+00
& 6.6e-01 & 1.4e+00
\\
\midrule
\multirow{2}{*}{Sphere}
& $E$
& 2.1e+02 & {\bf7.6e+01}
& 2.0e+02 & {\bf6.7e+01}
& 8.5e+02 & {\bf2.2e+02}
\\
& $\sigma$
& 2.0e+01 & 2.3e+01
& 1.6e+01 & 2.2e+01
& 7.3e+01 & 9.3e+01
\\
\midrule
\multirow{2}{*}{Squares}
& $E$
& 3.6e+04 & {\bf1.1e+04}
& 3.7e+04 & {\bf1.2e+04}
& 1.6e+05 & {\bf4.6e+04}
\\
& $\sigma$
& 1.6e+03 & 3.5e+03
& 4.6e+03 & 3.6e+03
& 1.9e+04 & 2.3e+04
\\
\midrule
\multirow{2}{*}{Trid}
& $E$
& 7.5e+08 & {\bf1.3e+08}
& 8.0e+08 & {\bf1.5e+08}
& 3.5e+09 & {\bf2.4e+08}
\\
& $\sigma$
& 8.9e+07 & 5.0e+07
& 4.5e+07 & 5.6e+07
& 4.9e+08 & 2.2e+08
\\
\midrule
\multirow{2}{*}{Trigonometric}
& $E$
& 1.1e+06 & {\bf1.8e+04}
& 9.3e+05 & {\bf4.1e+04}
& 1.7e+06 & {\bf1.1e+04}
\\
& $\sigma$
& 1.0e+05 & 2.5e+04
& 5.1e+04 & 3.8e+04
& 3.5e+05 & 1.7e+04
\\
\midrule
\multirow{2}{*}{Wavy}
& $E$
& 9.2e-01 & {\bf5.4e-01}
& {\bf3.1e-01} & 4.4e-01
& 9.6e-01 & {\bf5.9e-01}
\\
& $\sigma$
& 1.8e-02 & 6.5e-02
& 8.4e-02 & 3.0e-02
& 4.9e-02 & 9.0e-02
\\
\midrule
\multirow{2}{*}{Yang}
& $E$
& {\bf6.6e-11} & {\bf3.1e-32}
& {\bf9.0e-32} & {\bf3.1e-33}
& 1.7e+00 & {\bf3.4e-32}
\\
& $\sigma$
& 4.9e-02 & 9.0e-02
& 4.9e-02 & 9.0e-02
& 5.1e+00 & 4.4e-32
\\
\bottomrule
\end{tabular}
\end{table}

To rigorously validate the efficacy of TESALOCS, we conducted systematic comparisons across $20$ analytically defined benchmark functions with non-convex, multimodal landscapes~\cite{jamil2013literature}.
We compare the absolute error of the optimization result relative to the known exact global minimum for various existing optimization methods as part of our TESALOCS approach and without it (i.e., with random selection of the starting point).
For each considered optimization method we report the average value $E$ of the absolute error and its variance $\sigma$ for $10$ independent runs of each experiment.

All experiments were performed in a $100$-dimensional parameter space ($d = 100$) under strict computational constraints, limiting the number of objective function evaluations to $M = 10^4$ and with fixed rank $r=5$, sampling rate $N=2^{20}$, batch size $k=100$, elite size $k_{top}=10$.
The search space parameters $\bf{a}$ and $\bf{b}$ were selected according to benchmark function definitions.

In our initial experimental series, we evaluated widely-used gradient-based optimization methods~\cite{Dembo1983, Nocedal2006Numerical} from the SciPy library\footnote{
    See SciPy optimization module: \url{https://docs.scipy.org/doc/scipy/reference/optimize.html}
}: BFGS (Broyden–Fletcher–Goldfarb–Shanno algorithm); CG (Conjugate Gradient algorithm); L-BFGS-B (limited-memory BFGS variant with box constraints); Newton-CG; SLSQP (Sequential Least Squares Programming); and TNC (Truncated Newton Algorithm).
The computation results are reported in Table~\ref{tab:experiments_func1}, and Table~\ref{tab:experiments_func2}.
For each benchmark-method pair, bold formatting denotes the lowest absolute error.
Co-convergence is declared when both (random and our) approaches achieve errors $<10^{-8}$, then we consider both results to be the best.
For ease of comparison, table headers display victory counts, i.e. the total bold entries per methodology.
As follows from the tables, our approach leads to a better result more than $2$ times more often than the original method with random initialization.
At the same time, the variance for our method also turns out to be significantly lower than for the baseline in most cases.

Figure~\ref{fig:results_conv} delineates the optimization efficacy as a function of expended computational budget for the 100-dimensional Ackley, Rosenbrock, and Yang benchmarks as representative test cases.
The plotted trajectories reveal systematic improvements in solution precision when employing TESALOCS across all evaluated gradient-based algorithms (except for the L-BFGS-B method applied to the Rosenbrock function, which showed higher accuracy with the basic random initialization).

As noted above, our approach can be applied not only in conjunction with the gradient optimization method, but even with the gradient-free one.
To illustrate this point, we also consider well-known gradient-free methods~\cite{spall1992multivariate, kennedy1995particle} from the nevergrad library\footnote{
    See \url{https://github.com/facebookresearch/nevergrad}
}: NB (NoisyBandit method), PSO (Particle Swarm Optimization algorithm), and SPSA (Simultaneous Perturbation Stochastic Approximation algorithm).
We conducted completely similar experiments for these methods as for the gradient-based methods above and reported the corresponding results in Table~\ref{tab:experiments_func3}.
As it turns out, in this case the advantages of our method are even more significant: in almost all cases we observe an improvement relative to random initialization when using TESALOCS.
\section{Related work}
    \label{sec:related}
    Gradient-based optimization methods exhibit strong dependence on initial conditions due to their local search nature: the optimization trajectory becomes constrained by the basin of attraction surrounding the initial point and this is particularly problematic in high-dimensional landscapes where the probability of starting near global optima decreases exponentially with dimensionality~\cite{haji2021comparison, bian2024machine}.
Current approaches to address these limitations, which are based on the multi-start methods~\cite{povsik2012restarted, liu2024effective, sakka2025elegant}, can improve exploration but scale poorly with dimensionality due to exponential growth of required initial points.

These limitations of the gradient-based methods motivate our TESALOCS framework, which combines low-rank tensor sampling for global exploration with gradient-based local refinement.
By decoupling the global and local search phases through tensor decomposition, we overcome the dimensional scalability constraints while maintaining the precision of gradient methods, as shown in the previous experimental section.
Our approach directly addresses the initialization sensitivity through systematic tensor sampling and avoids local optima stagnation via discrete exploration of the parameter space.

The important innovation we propose is the use of a low-rank discrete tensor-based approach, i.e. TT-decomposition, as a sampler for the local method.
In recent years, the TT-decomposition~\cite{oseledets2011tensor} has gained significant attention due to the possibility of its successful use in multidimensional problems.
Recent advances in TT-based optimization have demonstrated remarkable success across various domains. Foundational work~\cite{sozykin2022ttopt} established TTOpt optimization algorithm via maximum volume principle, which is particularly effective in reinforcement learning scenarios.
Building on this,~\cite{batsheva2023protes} introduced PROTES -- a probabilistic sampling approach leveraging TT-decomposition for black-box optimization, achieving good results in high-dimensional discrete spaces.
Subsequent developments extended these methods to recommender systems~\cite{nikitin2022quantum}, adversarial attacks~\cite{chertkov2023tensortrain}, machine translation security~\cite{chertkov2024translate}, quantum circuits~\cite{selvanayagam2022global}, multiscale modeling~\cite{oseledets2016black}, robotics~\cite{shetty2022tensor}, neural activation analysis~\cite{pospelov2024fast}, etc., which demonstrates the versatility of the TT-based approach.

While these TT-based optimization methods excel in specific contexts, limitations remain: TTOpt, PROTES and other existing TT-approaches~\cite{dolgov2020approximation, novikov2021tensor, ryzhakov2022constructive, chertkov2023tensor, bharadwaj2024efficient} focuses solely on discrete spaces and lack integration with gradient-based refinement.
Our TESALOCS method addresses these gaps by combining adaptive TT-sampling with arbitrary gradient optimizers, enabling simultaneous handling of discrete and continuous variables while maintaining the dimensional efficiency of the TT-representation.
Thus, our hybrid approach extends the applicability of tensor methods while maintaining their computational advantages over traditional multivariate optimizers.

Combining local and global optimization strategies holds significant promise, as it merges the rapid convergence of gradient-based methods with the robustness of global algorithms in handling initial conditions.
In general, a combination of global and local (not necessarily gradient-based) methods is highly important in physics and material science~\cite{10.1063/1.3097197}, reinforcement learning~\cite{wang2022policy}, trajectory optimization~\cite{10605252}, drug therapy~\cite{Wu2023}, etc. This makes the proposed TESALOCS approach valuable for these areas.

\section{Conclusions}
    \label{sec:conclusions}
    We presented the algorithm TESALOCS for optimizing an unknown multidimensional function that combines the advantages of continuous methods and discrete optimization.
The key innovation lies in its two-stage architecture: discrete search for initial points and gradient-based (or even gradient-free) refinement.
Discrete optimization method is used to find a good starting point for any other method, and it ``sees'' only the best points (obtained after local optimization) and thus follows the maxima of the target function.

The potential scalability of our approach makes it adaptable for hyperparameter optimization in machine learning; for materials design and molecular dynamics tasks where both precision and speed are critical, etc.
Our future research will focus on the robustness analysis of the method, including under noisy objective function conditions, and its extension to the constrained optimization problems.
It is also of particular interest to consider other types of tensor networks instead of the TT-decomposition as the basis of the discrete optimization method used in our approach.

    \bibliographystyle{plainnat}
    \bibliography{refs}
\end{document}